\documentclass[12pt]{amsart}
\usepackage{epsfig,amssymb,epic,eepic,color}
\usepackage[all]{xy}
\usepackage{hyperref}
 \usepackage[T1]{fontenc}        
\newtheorem{thm}{Theorem}[section]

\newtheorem{prop}[thm]{Proposition}
\newtheorem{defn}[thm]{Definition}
\newtheorem{lemme}[thm]{Lemma}
\newtheorem{cor}[thm]{Corollary}

\newtheorem{remarques}[thm]{Remarks}

\newtheorem{rien}[thm]{}
 
\newcommand{\be}{\begin{enumerate}}
\newcommand{\ee}{\end{enumerate}}
\newcommand{\bi}{\begin{itemize}}
\newcommand{\ei}{\end{itemize}}

\def\R{\mathbb{R}}

\def\M{\mathbb{M}}

\def\om{\omega}

\def\ga{\gamma}

\def\al{\alpha}
\def\be{\beta}
 
\def\de{\delta}

\def\vp{\varphi}

\def\Si{\Sigma}

\def\ep{\varepsilon}

\def\nd{\noindent}
\def\bull{\hfill$\Box$\\}
\def\proof{\nd {\bf Proof.\ }}

\begin{document}
${}$
\vskip 1cm

\begin{center}
{\sc A proof of Reidemeister-Singer's theorem by Cerf's methods
\vskip 1cm 
 Fran\c cois Laudenbach
}

\today
\end{center}
\title{}
\address{Laboratoire de
Math\'ematiques Jean Leray,  UMR 6629 du
CNRS, Facult\'e des Sciences et Techniques,
Universit\'e de Nantes, 2, rue de la
Houssini\`ere, F-44322 Nantes cedex 3,
France.}
\email{francois.laudenbach@univ-nantes.fr}



\begin{abstract} Heegaard splittings and Heegaard diagrams of a
 closed 3-manifold $M$
are translated into the language of Morse functions with Morse-Smale
 pseudo-gradients
defined on $M$. We make use in a very simple setting of techniques
which Jean Cerf  developed
for solving a famous {\it pseudo-isotopy} problem. In passing, we show how to cancel the supernumerary
local extrema in a generic path of functions when $\dim M>2$. The main tool that we introduce 
 is  an {\it elementary
swallow tail lemma} which could be useful elsewhere.

\end{abstract}
\maketitle

\thispagestyle{empty}
\vskip .5cm
\section{Introduction}\label{intro}
\medskip
When speaking of Cerf's methods 
we refer to Cerf's  work 
in \cite{cerf} for the so-called 
{\it pseudo-isotopy}
problem. In a few words, the method consists of reducing  some 
isotopy problem to  a problem about real functions. It was created
 in the setting of high dimensional manifolds. However, some parts 
apply in dimension three as we are going to show. The purpose of this note is to present a proof 
of Reidemeister-Singer's theorem (as stated below) in this way.
I should say that Francis Bonahon, who like me  was educated 
in the Orsay Topology group of the seventies-eighties, wrote such a proof; 
but,  his notes are lost. The recent developments in Heegaard-Floer homology
drove me to make this proof available.
The concepts used in the next statement will be   explained in the course of  this introduction.
We always work in the $C^\infty$ category (also called the {\it smooth} category),
 for objects, maps and families of maps.
\begin{thm}{\bf (Reidemeister \cite{reidemeister}, Singer\cite{singer})}\label{r-s}
Let $M$ be a closed connected 3-manifold. 

\nd {\rm 1)} Two Heegaard splittings 
become isotopic after suitable stabilizations.
 
\nd {\rm 2)} More precisely, let $D_0, D_1$ be
 two Heegaard diagrams.  %
 Then there are stabilizations $D'_0,D'_1$
by adding pairs of cancelling handles of index 1 and 2,
such that one can pass from $D'_0$ to $D'_1$ by an ambient isotopy and  a finite sequence 
of  handle slides.\\
\end{thm} 

Strictly speaking, only the first item  is the statement of the Reidemeister-Singer theorem.
 A {\it Heegaard splitting } consists of a closed surface $\Si$ of genus $g$, called {\it Heegaard surface},
 dividing  $M$ into two handlebodies $H^-, H^+$. A {\it Heegaard diagram} is defined by
 more precise data, namely,
a handle decomposition of $M$ with:
\begin{itemize}
\item one 0-handle $B^-$ and $g$ handles of index 1 attached on the boundary   $\partial B^-$,
whose union forms $H^-$;
\item $g$ handles of index 2  attached on $\partial H^-$ and one 3-cell $B^+$, whose union forms 
$H^+$.
\end{itemize}
On the common boundary $\Si$ of $H^+$ and $H^-$, the Heegaard diagram specifies
 $g$ simple curves $\beta_1,...,\beta_g$ in $\Si$,
 mutually disjoint,  which are the cores of the attaching domains 
of the 2-handles; their complement in $\Si$ is a 2-sphere with $2g$ holes.
 It also specifies $g$ simple
curves $\al_1,...,\al_g$ which are the boundaries of the so-called {\it transverse} 2-cells\footnote{They are also called {\it compression discs}.}
of each 1-handle; 
the complement in $\Si$ of $\cup_j\al_j$ is also a 2-sphere with $2g$ holes. The other notions 
involved in Theorem \ref{r-s} will be only defined in the {\it functional  setting} considered below.\\

The statement  of Theorem \ref{r-s} can be 
  translated into the language of Morse functions
as follows. 
Recall that a  {\it Morse function}  $f$ is a smooth function whose critical points are non-degenerate; 
the famous {\it Morse lemma} states  that each critical point $p$ of $f$
belongs to a chart equipped with 
 so-called
{\it Morse coordinates}, meaning that $f-f(p)$ 
reduces  to  a quadratic form. 
Some non-classical facts  concerning the choice of these coordinates will be detailed in Section \ref{s3}.

A Morse function 
is said to be 
{\it ordered} if the order of the  critical values is finer than the order
of  their indices, namely $f(p)<f(p')$ whenever the index of the critical point $p$ is less than the index of 
$p'$. In dimension 3, an ordered Morse function gives rise to a Heegaard splitting 
by considering a level set whose level separates the index 1 and index 2 critical values.
Moreover, every Heegaard splitting is obtained this way.
Along a path of ordered Morse functions the  Heegaard surface moves by isotopy.

Stabilizing a Heegaard splitting 
consists of  creating a pair of critical points of index 1 and 2 at a level
keeping the ordering.
Thus, 
item 1 of Theorem \ref{r-s} is a consequence of  Theorem \ref{order-path}, for which it is necessary 
to speak of {\it genericity}.

\begin{rien} {\bf Genericity I.} \label{genI} {\rm Given two Morse functions $f_0,f_1: M\to\R$,
the following property is {\it generic} (in Baire's sense) for the paths of functions $\left(f_t\right)_{t\in [0,1]}$ joining them:
\begin{itemize}
\item for all $t\in[0,1]$ apart from finitely many  {\it exceptional} values 
$t_j$,  the function  $f_t$ is Morse;
\item 
 for $\de >0$ small enough, 
$f_{t_j+\de}$ has one more or one less pair of critical points  than $f_{t_j-\de}$; in the first 
(resp. second) case, 
$t_j$ is  called a {\it birth time} (resp. a {\it cancellation time});
\item  the critical points of $f_{t_j}$ are all 
 non-degenerate  except  one whose Hessian has corank  1; this point 
 will be said a {\it cubic critical point}.
\end{itemize}
For short, when speaking of a {\it generic} path of functions, it will be understood a path as above. 

 In this note, all  genericity argument  follow
from Thom's {\it  transversality theorem in jet spaces} as it is  in his article
on singularities \cite{thom56} (see also  \cite{hirsch}, or \cite{lauden} 
where the generic paths of real functions are  explicitly considered). 
In Section \ref{s2}, we shall specify which transversality is involved in the above genericity of paths.\\

The next theorem  is mainly due to Jean Cerf (\cite{cerf}, chap. V \S I)\footnote{Strictly speaking, only the first sentence is stated in Cerf's article. The complement  follows from his lemma about the {\it uniqueness of births}
(valid in dimension greater than 1 only).}. 
}
\end{rien}

\begin{thm}\label{order-path} Let $M$ be a closed connected manifold of any dimension $n$. 
Given two ordered Morse functions $f_0,f_1$ on $M$, they are joined 
by a generic path of functions 
$\left( f_t\right)_{t\in[0,1]}$ 
such that, for every $t\in[0,1]$ 
outside of   a finite set $J=\{t_1,\ldots, t_q, t_{q+1},\ldots, t_{q+q'}\}$, 
$f_t$ is an ordered Morse function. Moreover, $t_1, \ldots, t_q$ are 
birth times and lie in $\left(0,  \frac13\right)$;  and  $t_{q+1},\ldots, t_{q+q'}$
are cancellation (or death) times and lie in $\left(\frac 23, 1\right)$. 

In particular in dimension 3,  a level set of $f_{1/2}$ whose level separates
the  index 1 and index 2 critical values is a Heegaard splitting that is a common
stabilization, up to isotopy, of those associated with $f_0$ and $f_1$. \\
\end{thm}

We now turn to the second part of Theorem \ref{r-s}.
In order to speak of handle decomposition and handle sliding,  it is useful to 
consider a Morse function $f$ equipped with a {\it pseudo-gradient}. 

\begin{defn} \label{pseudo} Given a Morse function $f$,  a smooth vector field $X$ on $M$ is said to be a (descending)
pseudo-gradient for $f$ if the two following conditions hold:
\begin{itemize}
\item the Lyapunov inequality\footnote{This sign convention   is used for instance by R. Bott p. 341 in \cite{bott}. } $X\cdot f<0$ away from the critical locus;
\item  at each
critical point $p$ the Hessian 
of $X\cdot f$ is negative definite (notice that $X\cdot f\leq 0$
everywhere). 
\end{itemize}
\end{defn}

Local data of pseudo-gradients generate a global pseudo-gradient by using a partition of unity.
It is easily checked that the zeroes of $X$ coincide with the critical points of $f$
and are {\it hyperbolic}\footnote{That is, if $p$ is a zero of $X$ the eigenvalues of the linearized vector 
field at $p$ have a non-zero real part.}. Thus, according to the {\it stable/unstable manifold theorem} 
(see \cite{brin}), with each zero $p$ of $X$
 there are  associated stable and unstable manifolds, also called {\it invariant manifolds} and 
denoted respectively  by $W^s(p,X)$ and $W^u(p,X)$. A point $x\in M$ belongs to 
$W^s(p, X)$ if $X^t(x)$ tends to $p$ as $t$ tends to $+\infty$; here, $X^t$ denotes the flow of $X$.

The unstable manifold is diffeomorphic to $\R^i$,
where $i$ is the index of $f$ at 
$p$, 
and the stable manifold is diffeomorphic to $\R^{n-i}$; moreover, 
$p$ is a non-degenerate maximum (resp. minimum) of the restriction of $f$ to $W^u(p,X)$
(resp. $W^s(p,X)$).

 Given the Morse function 
$f$,  Smale \cite{smale} proved   that, generically,  
   all  invariant manifolds of a pseudo-gradient of $f$ are mutually transverse\footnote{An ordered Morse function
$f$ with a Morse-Smale pseudo-gradient $X$ gives rise easily to a handle decomposition.}. Today, 
such a pseudo-gradient is said to be {\it Morse-Smale}.

According to Whitney \cite{whitney}, if $p$ is a cubic critical point of $f$, there are coordinates
$(x,y)\in \R\times\R^{n-1}$, which we call
{\it Whitney coordinates}, where $f$ reads:
$$f(x,y)= f(p)+ x^3+q(y).
$$
Here, $q$ is a non-degenerate quadratic form on $\R^{n-1}$. For  a reason
which will be explained in  
 \ref{fam-pseudo}, we require a pseudo-gradient $X$ for $f$ to coincide 
with $-\nabla_gf$ near the cubic critical point $p$,
where $g$ is the Euclidean metric of one system of Whitney coordinates.

Given a generic path of functions $\left( f_t\right),\, t\in [0,1],$ it can be enriched with a smooth path of vector fields $\left (X_t\right)$, such that $X_t$ is  a pseudo-gradient of $f_t$ for all $t\in [0,1]$.

\begin{rien}{\bf Genericity II.} \label{genII}{\rm The following property is generic  for the paths of pairs
$\left(f_t,X_t\right)_{t\in [0,1]}$:
\begin{itemize}
\item the path of functions is generic in the sense of \ref{genI};
\item for every $t$,  there is no $X_t$-orbit from a critical point index $j$ of $f_t$ to a critical point
 index $i$ if $j<i$ (briefly said: no $j/i$ {\it connecting orbit}  if $j<i$);
 \item for every $t$ outside of a finite set $K=\{t_1,\ldots,t_r\}\subset(0,1)$ of Morse times\footnote{A cubic
  point of index $i$ could be connected to a Morse point of index $i$ at a lower level.},
  there is no $i/i$
 connecting orbit of $X_t$;
 \item for each $t_k\in K$, exactly one orbit $\ell_k$ of $X_{t_k}$
 connects two critical points $p$ and $p'$ 
 having the same index; moreover, for each $x\in \ell_k$, we have:
$$T_x\ell_k=T_xW^u(p,X_{t_k}) \cap T_xW^s(p',X_{t_k})\, ,$$
and $t\mapsto X_t$ crosses transversely  at time $t_k$ the codimension-one stratum
 of the space of pseudo-gradients having a connecting orbit between two critical points with the same index.
\end{itemize} }
\end{rien}
For short, such a   path $\left(f_t,X_t\right)_{t\in [0,1]}$ is said to be {\it generic}. 
For $t_k\in K$,
one says that a {\it handle sliding} happens at time $t_k$. The effect of 
a handle sliding  on the so-called 
{\it Morse complex} is described by  J. Milnor
(see Theorem 7.6 in \cite{h-cob}). 

The argument for genericity  in \ref{genII} is elementary once the first item is assumed. It  relies on the classical  transversality theorem
applied to a  $(j-1)$-sphere moving with $t$ with respect 
to a fixed $(n-i-1)$-sphere,  $j\leq i$, in an $(n-1)$-dimensional manifold.

Now, the statement of item 2) in Theorem \ref{r-s} can be  translated into the next one. Following
  M. Morse  \cite{morse},   a function with only two local extrema is be
  said to be {\it polar}.

\begin{thm}\label{ord-f} Let $M$ be a closed connected manifold of dimension\footnote{The statement
also holds in dimension 2 with a different proof (see \cite{kudry}, \S 8).  It is obvious in dimension 1.} $n>2$.
Given two ordered polar  Morse functions $f_0,f_1$ equipped with respective Morse-Smale 
pseudo-gradients $X_0,X_1$, there exists 
a generic path of pairs $(f_t,X_t)_{t\in[0,1]}$, where the vector field 
$X_t $ is a pseudo-gradient for the function
$f_t$, so that the following holds:
for every $t\in [0,1]$ outside of  a finite set, 
 $f_t$ is an ordered  polar Morse function and $X_t$ has no $i/i$ connecting orbit. 
 The excluded values of $t$ are the times of birth first, then handle sliding and finally cancellation.
\end{thm}

A direct  proof of Theorem \ref{order-path} is given in Section \ref{s2} without any reference to Cerf's work.
It mainly follows  from  Lemma \ref{reorder} which offers an efficient process
for crossing critical values. 
The proof of Theorem  
 \ref{ord-f} will be given in Section \ref{s4} and  uses a few
  technical lemmas, including the {\it elementary swallow tail lemma} and the {\it elementary lips lemma}.
  Since they could be useful in a more general setting, they are written with index assumptions which are 
  more 
  general than necessary here. These lemmas are proved in Section \ref{s3}.
  
  \section{Proof of Theorem \ref{order-path}} \label{s2}
  
  The main tool will be the next lemma.

\begin{lemme} {\bf (Decrease of a critical value)}\label{reorder}
 Let $f:M\to \R$ be a Morse function, let $X$ be a pseudo-gradient
for $f$ and let $p$ be  a critical point of index $k$.
Assume  that the unstable manifold $W^u(p,X)$ contains a closed smooth $k$-disc
$D$ 
whose boundary lies in a level set $f^{-1}(a),\ a<f(p)$. 
Then, 
for every $\ep>0$ with $a+\ep<f(p)$,
 there exists a path  $(f_t)_{t\in[0,1]}$ of Morse functions such that
$f_0=f$, $f_1(p)=a+\ep$ and $X$ is a pseudo-gradient of $f_t$ for every 
$t\in[0,1]$. Moreover, the support of the deformation may be contained 
in an arbitrarily small neighborhood $W$ of  $D$ in $M$. 
\end{lemme}
 
 Note that, when $k=0$, $W^u(p,X)$  has an empty intersection 
 with the open sub-level set  $f^{-1}\bigl((-\infty, f(p))\bigr)$. So, the condition of the lemma is fulfilled
 and the conclusion  allows us to decrease arbitrarily the value of a local minimum.

The  lemma above holds true, with the same proof, in  a family whose data 
$(f,p,D,a)$ depend smoothly on a parameter $s\in \R^m$
and fulfill  the same assumptions for every $s$. Moreover, 
$f$  only has  to be a Morse function in a neighborhood of $D$. 
In particular, it applies  to non-generic  functions or pseudo-gradients.\\

\nd  {\bf Proof.} The case where $p$ has index 0 is left to the reader. 
Hereafter, assume $k>0$. Set $n=\dim M $ and $c=f(p)$. For $\eta$ small enough, there exists a 
closed $(n-k)$-disc $D'$ in the stable manifold $W^s(p,X)$, with $D'\subset W$,
 whose boundary lies in $f^{-1}(c+\eta)$.
  Let $U$ be a tubular neighborhood of radius $\de$
of $\partial D$ in $f^{-1}(a)$.  For $\de$ small enough with respect
to $\eta$, every half-orbit of $X$ 
ending in $U$ is contained in $D$ or crosses $f^{-1}(c+\eta)$. Define $\mathcal M $ as 
the union of $D$, $D'$
 and all segments of $X$-orbits starting from points in $f^{-1}(c+\eta)$ and ending  in $U$; for a small 
 $\de $, we have $\mathcal M\subset W$. 
Its boundary 
is made of three parts, two horizontal parts 
$\mathcal M\cap f^{-1}(a)$
and $\mathcal M\cap f^{-1}(c+\eta)$, and the lateral boundary 
$\partial_\ell\mathcal M$ which is tangent to $X$.  There are two corners 
in the boundary of $\mathcal M$, each  being diffeomorphic to a product of spheres 
$S^{k-1}\times S^{n-k-1}$ (where $k=index(p)$); one is the boundary of $U$, trivialized as the sphere normal bundle 
$\partial U\to \partial D$; the other corner is $\partial_\ell\mathcal M \cap f^{-1}(c+\eta)$
and is diffeomorphic to the first one by the flow of $X$.

Let $N$ be a small
collar neighborhood of $\partial_\ell\mathcal M$ in $\mathcal M$; it is diffeomorphic to a product
$$N\cong S^{k-1}\times S^{n-k-1}\times[0,1]\times[a,c+\eta].
$$
If $(x,y)$ are the coordinates of $R:=[0,1]\times[a,c+\eta]$, 
the product structure of $N$ is chosen so that the level sets of $f$ in $N$ are $\{y=const.\}$
and the vertical lines directed by $\partial_y$
 are  tangent 
to the orbits of $X$.

For constructing $f_1$ we keep $f_1=f$ outside of $\mathcal M $ and change the level set foliation
as said below. The level set foliation of $f_1$ coincides with the one of $f$ 
in the complement of $N$ in $\mathcal M$. Inside  $N$,
it is obtained by  replacing the horizontal foliation of $N$ with a new one
which is still transverse to the vertical lines, 
 is still horizontal near the boundary,
and puts $f^{-1}(a+\ep)\cap \{x=0\}$ on the same leaf as $f^{-1}(c)\cap \{x=1\}$. The new foliation in $N$
 is the pullback of a foliation of $R$ by the standard projection  
(see figure 1). 
The value of $f_1$ is now well-defined. 
\begin{center}
\hskip 0cm \includegraphics[scale=0.6]{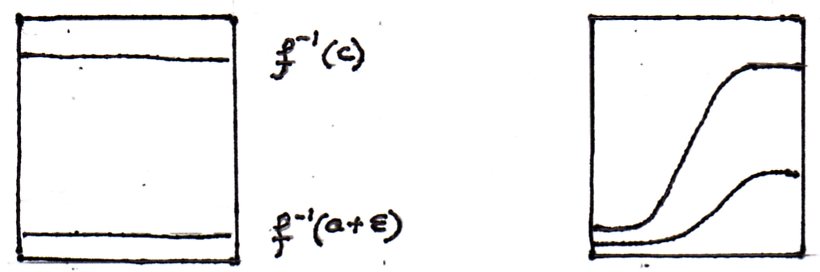}
\centerline{Figure 1A \hskip 4cm Figure 1B}
\end{center}
Moreover, it is easy to interpolate
this construction for $t$ varying in $[0,1]$. \bull\\

\begin{cor} \label{reordering_0}Let $(f_0,X)$ be a Morse function with a pseudo-gradient
having no $j/i$ connecting orbit, $j<i$. 
Then there exists a path $(f_t)_{t\in[0,1]}$ of Morse functions issued from
$f_0$  such that $f_1$
is ordered and the same vector field $X$ is a pseudo-gradient of $f_t$ for every $t\in [0,1]$.
\end{cor}

\medskip

\proof 
If the function is not ordered,
there is a pair of critical points $(p, q)$ with $index(p)<index(q)$ 
and $f(p)\geq f(q)$.
Choose such a pair  so that $f(p)$  
is  minimal among all similar 
{\it unordered pairs}. By this choice every orbit of  $W^u(p, X)$ crosses a level set 
below $f(q)$; if not, one of them  ends at a critical point $p'$. By assumption on $X$ we have
$index(p')\leq index(p)<index(q)$
and  $f(p)>f(p')\geq f(q)$,
 contradicting the assumption on the pair 
$(p,q)$. Then, lemma \ref{reorder} applies and yields a new Morse function 
which has the same pseudo-gradient $X$ and  at least one unordered  pair less than $f$.
Arguing this way recursively, the corollary is proved.
\bull

Before proving Theorem \ref{order-path}, it is useful to specify which transversality is involved
in a generic path in the sense of \ref{genI} and what a {\it birth path} is. A path of functions 
$\left(f_t\right)$ may be thought 
of as a smooth function $F:[0,1]\times M\to \R, \ (t,x)\mapsto f_t(x)$. We now consider  the 
$r$-jet spaces $J^r([0,1]\times M,\R)$ for $r=1,2$ and their submanifolds
$\Si^1$ and $\Si^{1,1}$ defined as follows (here, we  are using  the so-called Thom-Boardman
notation). The first one, $\Si^1$, is made of the 1-jets $(a,j^1g)$ where  $a\in [0,1]\times M$ and $g$ is
a germ at $a$ of function $(t,x)\mapsto g(t,x)$ such that $\partial_xg (a)=0$. The second one, $\Si^{1,1}$,
is made of the 2-jets $(a,j^2g)$ such that:
\begin{itemize}
\item $dg_xg(a)=0$ and $j^1g$ meets $\Si^1$ transversely; 
\item the (germ of) curve $\left(j^1g\right)^{-1}(\Si^1)$ passes through $a$  and is tangent
 to the kernel of $\partial_xg(a)$, which is the factor $\{t=t(a)\}$.
\end{itemize}

According to Thom \cite{thom56}, generically $j^1F$ is transverse to $\Si^1$ and $j^2F$ is transverse 
to $\Si^{1,1}$. Thus,  the critical locus of $f_t$ when $t$ runs in 
$[0,1]$, which is   $\left(j^1F\right)^{-1}(\Si^1)$, is a   smooth curve; 
and the isolated points $\left(j^2F\right)^{-1}(\Si^{1,1})$ are the cubic critical points.
By making a diffeomorphism  $C^\infty$-close to $Id$ act on  $[0,1]\times M$,
 it is possible to
move the cubic critical points so that their $t$-coordinates are distinct. In particular, 
the properties in \ref{genI} hold true generically.

Moreover, if $(t_0,x_0)$ is a cubic critical point, thanks to the information on the 3-jet\footnote{The transversality of $j^2F$  to $\Si^{1,1}$ at $(t_0,x_0)$ is an open  condition on the 3-jet.}
of $F$ at $(t_0,x_0)$, it is possible to write a normal form of $F$ on a neighborhood of $(t_0,x_0)$.
This follows easily from the normal form of {\it cusps} established  by H. Whitney in \cite{whitney}
for generic maps from plane to plane. Precisely, there are adapted coordinates $(t,x)=(t,y,z)$, with 
$y\in \R^{n-1}, z\in \R$,
which we call {\it Whitney coordinates},
where $F$ reads:
$$\quad F(t,x)= F(t_0,x_0) + z^3\pm (t-t_0)z +q(y)
$$
 Here, $q$ is a non-degenerate quadratic form on $\R^{n-1}$, $\pm=-$ if $t_0$ is a birth time and  $\pm= +$
 if $t_0$ is a cancellation time. If $t_0$ is a birth time,  we immediately derive from the model
 that, for $\de>0$ small enough, the given generic path of functions, restricted to
 $[t_0-\de,t_0+\de]$, is a birth path in the following sense.

\begin{defn} \label{birth} A birth path is a generic path of functions $\left(f_t\right)_{t\in [t_0-\de,t_0+\de]}$
such that there exists a path of 
cylinders $B_t\cong D^{n-1}\times [-1,+1]$ embedded in $M$
with the following properties for every $t\in [t_0-\de,t_0+\de]$:
\begin{itemize}
\item $D^{n-1}\times \{\pm 1\}$ (the top and bottom of $B_t$) lie in two level sets of $f_t$;
\item  the restriction of $f_t$ to $\partial D^{n-1}\times [-1,+1]$ has no critical points;
\item $f_t \vert B_t$ is semi-conjugate to  the function $c_{t_0}^t(y,z):= z^3-(t-t_0)z +q(y)$. 
\end{itemize}
Here, a semi-conjugation stands  for an embedding $\vp_t: B_t\to\R^n$, 
depending smoothly on $t$, covering the origin of $\R^n$ and such that 
$c_{t_0}^t\circ \vp_t=f_t\vert B_t$ up to a rescaling of the values. 
The index of $q$ is  called the {index of the birth}. 
\end{defn}
The function $f_{t}$ has no critical points in $B_{t}$ when $t_0-\de\leq t<t_0$
whereas, for $t_0<t\leq t_0-\de$,
 $f_{t}$ has a pair of critical points in $B_{t}$ of respective index $i, i+1$ if $i$ is the index of the birth. 

\begin{remarques}\label{birth-uniq}{\rm 
 1)  If $f_0$ is a Morse function given with a cylinder $B_0$ on which $f_0$ induces the height function, then
$f_0$ is the beginning  of a birth path with $t\in [0,2\de]$ which is supported in $B_0$ in the sense  that 
 the path is stationary
outside of $B_0$. Indeed,  $f_0\vert B_0$ is semi-conjugate to any  function without critical point,
for instance $(y,z)\mapsto z^3+ \de z + q(y)$; thus, it is allowed to plug
 the functions $c_\de^t, \ t\in [0,2\de]$,
by taking a suitable semi-conjugation $\vp_t: B_0\to\R^n$. This birth path is said to be  {\it elementary} 
(compare with a similar definition  in  Cerf \cite{cerf} chap. III).

2) Any birth path issued from $f_0$ associated with a path of cylinders $\left(B_t\right)_{t\in [0,2\de]}$
is homotopic  to an elementary birth path among the birth paths starting from  $f_0$. This is done by using
an extension of the isotopy $B_0\to B_t$.
}
\end{remarques}


\begin{lemme}{\bf (Shift of birth)}\label{shift} ${}$

{\rm 1)} Every generic path of functions on $M$ is homotopic relative to its end points to a generic path
where the birth times appear before the cancellation times. More precisely, the following holds.

{\rm 2)} Let $\left(h_s\right)_{s\in [0,1]}$ be a generic path of functions which are Morse for all time except
one cancellation time.
Let $\left(\beta ^1_t\right)_{t\in[0,2\de]}$ be a birth path starting from the Morse function
$h_1$ with associated cylinders
$\left(B^1_t\right)_{t\in[0,2\de]}$. Then  there is a smooth family, parametrized by $s\in [0,1]$,
 of birth paths $\left(\beta^s_t\right)_{t\in[0,2\de]}$, starting  from $h_s$ 
  with  associated cylinders $\left(B^s_t\right)_{t\in[0,2\de]}$ which coincide with the given cylinders when 
 $s= 1$. 

Moreover, if 
  $\dim M>1$,  the same holds  true for  any generic path $\left(h_s\right)_{s\in [0,1]}$.
Moreover, it is possible to choose the cylinders $B^0_t$
as neighborhoods of any given regular 
point of $h_0$.
\end{lemme}


\nd {\bf Proof of 2)$\Rightarrow$1).} The composed  path $\left(h_s\right )_{s\in[0,1]}*\left(\beta^1_t\right)_{t\in[0,2\de]}$ is homotopic, relative to its end points, to the composed 
path  $\left(\beta^0_t\right)_{t\in[0,2\de]}*\left(\beta_{2\de}^s\right )_{s\in[0,1]}$. In general, this 
composition is only piecewise smooth at the gluing point.
But we are free to modify the  parametrization of the composed path; if the two paths
entering the composition are stationary near their common end point, then the composed path is smooth.

The new path from $h_0$ to $\beta^1_{2\de}$ has one  birth time appearing before one cancellation time.
By arguing this way recursively  one proves 1).\\

\nd {\bf Proof of 2).}
Given the cylinder $B^1_0$, one chooses a smooth family of cylinders $\left(B_0^s\right)_{s\in [0,1]}$ in $M$ ending to $B^1_0$ 
and so that 
$h_s$ induces the standard horizontal foliation $D^{n-1}\times\{pt\}$ of $B^s_0$ for every $s\in [0,1]$. 
 This is possible in any positive dimension since we are free to move $B^s_0$ away from the critical set 
 of $h_s$, even at the cancellation time.
 Thanks to an extension of isotopies,
  we get 
a 2-parameter family of diffeomorphisms $\psi_t^s: B_0^s\to B^1_t$, $s\in [0,1], t\in[0,2\de]$,
preserving the horizontal foliation near the boundary and
such that $\psi^1_0= Id$. Then, define 
$$
\beta_t^s=\left\{
\begin{array}{l} 
h_s {\rm \ outside\ of\ }B_0^s\\
\beta^1_t\circ\psi_t^s\ {\rm in\ } B_0^s, \ {\rm up\ to\ some\ rescaling.}
\end{array}
\right.
$$ 
The rescaling is needed for making the two definitions match along the boundary of $B_0^s$.
When $t=1$, this is an elementary birth path issued from $h_1$. According to  Remark \ref{birth-uniq} 2),
it is homotopic to $\left(\beta_t^1\right)$ relative to $h_1$. This proves the first part of 2).
In case $\dim M>1$,  the critical locus is non-separating and the last statement of 2) follows.\bull\\

\begin{rien}{\rm {\bf Proof of Theorem \ref{order-path}.} The case $\dim M=1$
is left to the reader. Hereafter, $\dim M$ is assumed to be greater than 1. Given two ordered Morse functions
$f_0,f_1$, there exists a generic path $(f_t)_{t\in[0,1]}$ where $f_t$ is Morse 
for every $t\in [0,1]$ outside of a finite set $J$. 
Decompose $J= J_+\cup J_-$ where $J_\pm$ is the set of birth/cancellation times and apply
Lemma \ref{shift}.  The birth times $J_+$ can be shifted to the left, say in $[0,t_0]$,
and the cylinders
of birth can be located at the right level according to the index of the birth
so that all Morse functions in $[0,t_0]$ are ordered.
Similarly,  the cancellation times can be shifted to the right, say in $[t_1,1]$, and the cancellation cylinders can be chosen so that all Morse functions in $[t_1,1]$ are ordered.
Thus, $f_t$ is a Morse function for every $t\in [t_0,t_1]$ and is ordered for $t=t_0,t_1$.

  Choose  pseudo-gradients  $X_t$ for $f_t$.
 We may assume   $(X_t)_{t\in [t_0,t_1]}$   
 in the sense of \ref{genII}. Thus,  the pseudo-gradient $X_t$ has no $j/i$ connecting orbit with $j\leq i$
for all $t\in [t_0,t_1]$ outside of  a finite set $K\subset (t_0,t_1)$ (times of $i/i$ connecting orbits).

Apply corollary \ref{reordering_0} to the functions $f_{t_k}, \ t_k\in K$, and deform the path of functions accordingly, that is: keep the same path $\left(X_t\right)$ as path of pseudo-gradients
and ask the deformation to be 
 stationary on the complement of small neighborhoods of the $t_k$'s.
 After that deformation,
the functions $f_{t_k}, \ t_k\in K,$ are ordered and, for every  $t\in (t_k,t_{k+1})$,
the vector fields $X_t$ is 
has no   $j/i$ connecting orbit with $j\leq i$. This   also holds true on the intervals
$(t_0,\inf K)$ and $(\sup K, t_1)$
on the left and right  of $K$. So, we are
reduced to reorder a path of  Morse functions  equipped with pseudo-gradients which 
have no $j/i$ connecting orbits, $j\leq i$,  for every time. The reordering is then 
 obtained by applying the one-parameter version of 
 Lemma \ref{reorder}. This finishes the proof of item 1) in Theorem \ref{r-s}.\bull
}\end{rien}

\section{The elementary swallow tail lemma and similar results} \label{s3} 

Before proving Theorem \ref{ord-f} and, hence, item 2) in Theorem \ref{r-s},
 we need to state 
 some lemmas:
  first, a very particular case  of  the {\it swallow tail lemma}\,;
  next,  a very particular case of the {\it lips} lemma (or
 {\it uniqueness of death} according to \cite{cerf});
 finally, the {\it cancellation theorem}\footnote{Also referred simply as the {\it cancellation lemma}.}
  of Morse \cite{morse}
(see also J. Milnor \cite{h-cob}, Section 5). 
  
 We state them 
 by means of Cerf graphics. Recall that the {\it Cerf graphic}
  of a path of functions $\left(f_t\right)_{t}$ is the part of $\R^2$
 whose intersection with $\{t\}\times \R$ is the set of critical values of $f_t$.
 
The three proofs are very similar, by reduction to the one-dimensional case 
where they become easy.
Only the proof of the elementary swallow tail lemma is detailed here since the three proofs can be performed 
in the same way\footnote{Such a proof of Morse's cancellation theorem is now available in \cite{cancellation}.}. 

We begin with  useful conjugation lemmas. The first one is likely well-known, the next ones could be 
less 
 classical.

\begin{lemme} \label{conjugation}
Let $V$ be a manifold and $V'$ be a compact submanifold. Two germs of smooth functions $f$ and $g$ 
along $V'$ whose restrictions  to $V'$ coincide and have no critical points are isotopic relative to $V'$. 
Moreover, if $f = g$ near a compact set $K \subset V'$, the isotopy may be stationary near $K$ in $V$.
  This statement holds true with parameters in a compact set.
\end{lemme}

\nd{\bf Proof.} The path method of J. Moser \cite{moser} is available; it is explained below in our setting. 
Look at the path of germs 
 $t\in [0,1]\mapsto f_t:= (1-t)f+t g$ and search for an isotopy $\left( \vp_t\right)_{t\in[0,1]}$
 of $V$, with $\vp_0=Id$,
 satisfying the conjugation equation of germs along $V'$:
 $$(1)\quad\quad \begin{array}{l}
 f_t\circ\vp_t= f,\\
 \vp_t(x)=x \ {\rm for \ every\ } x\in V'.
 \end{array}
 $$
 The infinitesimal generator $Z_t$ has to satisfy the derived equation:
 $$(2)\quad\quad \begin{array}{l}
 df_t(x)\cdot Z_t(x)+g(x)-f(x)=0,\\
 Z_t(x)=0 \ {\rm for \ every\ } x\in V'.
 \end{array}
 $$Conversely, if $Z_t$ is a time depending vector field which is a solution of (2) near $V'$,
  its ``flow'' is defined until $t=1$  on a small neighborhood of $V'$ and solves the conjugation problem. 
 
 Here is a solution of Equation (2) by using an auxiliary Riemannian metric:
 $$Z_t= (f-g)\frac{\nabla f_t}{\vert \nabla f_t\vert^2}.
 $$
 The same proof holds for the relative statement and with parameters. \bull\\

\begin{lemme}{\bf (The $\mathfrak{ M J}^2$ lemma.)}\footnote{We learnt this proof of Morse's lemma 
from J. Mather on the occasion of a lecture in  Thom's seminar at IH\'ES (Bures-sur-Yvette), Dec. 1969.}
\label{mj2}
Let $\mathfrak F$ be the ring of germs of smooth functions at $0\in \R^n$ and let $\mathfrak M$ be its
unique  
maximal ideal of germs vanishing at 0. Given $f\in \mathfrak M$, its Jacobian ideal is the ideal $\mathfrak J=\mathfrak J(f)$ generated by 
the first partial derivatives of $f$. Consider a germ $h$ in the product ideal $\mathfrak{ M J}^2$. 
Then there is a $C^\infty $ diffeomorphism $\vp$ such that 
$(f+h)\circ \vp= f$.
\end{lemme}
For instance, 
take a  germ $f$ of Morse function with $f(0)=0$; it  reads
$f= q+r$  where $q$  is a non-degenerate quadratic form and  $r$  belongs to $\mathfrak M^3$. Since 
$\mathfrak J(q)=\mathfrak M$, 
the lemma implies that $f$ is conjugate to $q$, which 
is exactly the  statement of 
Morse's lemma.\\

\nd {\bf Sketch of proof.}\footnote{A detailed proof may be found in \cite{coursx}.} As in Lemma
\ref{conjugation}, we use the path method. Setting $f_t= f+th$, one searches for a family 
of local diffeomorphisms $\vp_t$, $t\in [0,1]$, such that $f_t\circ \vp_t= f$. This amounts
to find local vector fields $Z_t$ vanishing at the origin such that $df_t(x)\cdot Z_t(x)+ h(x)=0$; 
this consists of decomposing $h$ in the Jacobian ideal $\mathfrak J_t$ of $f_t$ with coefficients in 
$\mathfrak M$. The main point is that 
$\mathfrak J_t=\mathfrak J_0$ for all $t$. Indeed, $\left(\frac{\partial f_t}{\partial x_i}\right)= 
A_t\left(\frac{\partial f_0}{\partial x_j}\right)$ where the matrix $A_t$  equals the Identity matrix  
 modulo $\mathfrak M$.
Thus, $A_t$ is invertible, and  
 a decomposition of $h$ in $\mathfrak J_0$ with coefficients in $\mathfrak M$ yields the 
 wanted decomposition. 
${}$\bull

The same proof works with parameters $s\in \R^m $ and in a relative form: {\it Let
 $\left(f^s\right)_{s\in D^m}$ be 
 a family, parametrized by the $m$-ball, of  germs of Morse functions $(\R^n,0)\to\R$ 
whose Hessians at $0$ are denoted by $q^s$. Assume $f^s=q^s$
for every $s\in \partial D^m$. Then there is a family of local diffeomorphisms  $\vp^s$ such that
$f^s\circ \vp^s=q^s$ and $\vp^s= Id$ when $s\in \partial D^m$.}

In the same setting, if $f$ is given a local unstable manifold $W: =W^u(0, X)$, a system of Morse coordinates $x=(y,z)$ are said to be {\it adapted} to $(f,W)$ if  $f(x)=-\vert y\vert^2+\vert z\vert^2$
and $W=\{z=0\}$.

\begin{cor}\label{unstable-adapted} Given such data $f$ and $W$ the following holds.

1) There exist Morse coordinates adapted to $(f,W)$. (This claim also holds  with parameters.)

2) Two such systems of Morse coordinates can be joined, up to a permutation of the coordinates
by a one-parameter family of adapted Morse coordinates\footnote{We are hiding some acyclicity here
(compare \cite{chenciner}); but, the space of Morse coordinates is not acyclic, due to the isometry group $O(i,n-i)$. }.
\end{cor}

\nd{\bf Proof.} 1) The restriction  of $f$
to $W$ has a non-degenerate maximum. By Morse's lemma we have Morse coordinates $y$ of $W$ 
so that $f(x)= -\vert y\vert^2$ if $x\in W$. Complete the coordinates $y$ to  local coordinates $(y,z')$ of 
$(\R^n,0)$ so that $W=\{z'=0\}$ and the $z'$-space  is the orthogonal of the 
$y$-space with respect to $d^2f(0)$. Let $\left(f^y\right)$ be
 the family  of the restrictions of $f$ to the slice $\{y=cst\}$. For $y=0$, the 
function $f^0$ is Morse and its  critical point is $z'=0$. By the implicit function theorem, there is a smooth 
map $y\mapsto z'=k(y)$ such that $f^y$ is Morse with critical point at $k(y)$ (for every $y$ close to $0$).
Apply  the change of variables $(y,z)= (y, z'-k(y))$ so that  the critical point of $f^y$ becomes
 $z=0$ for every
$y$. By  a linear transformation in each slice, we may assume the Hessian of $f^y$ to  be
constantly  equal to $\vert z\vert^2$. The wanted Morse coordinates are now 
given by applying Morse's lemma 
with parameters to the family $\left(f^y\right)$.

2) We first connect the two given Morse coordinates by a path of coordinates which are only adapted to
$W$. Then, this path is modified by applying Morse's lemma  with parameters in the relative form.

\begin{rien}{\bf  Pseudo-gradients for birth path.}  \label{fam-pseudo} {\rm
To avoid raising  some problems in bifurcation theory of vector fields 
we adopt a still more restrictive
definition of pseudo-gradients\footnote{We could ask the path 
$\left(X_t\right)$ to present a bifurcation of type
{\it saddle-node}  along a birth/cancellation path.} 
than in \ref{pseudo}. 
This is allowed since 
we are free to choose our pseudo-gradients.

Recall from 
\ref{birth} (with slightly different notation) that a birth path at time $t_0$ consists of 
a generic path of functions $\left(f_t\right)_{t\in (t_0-\de,t_0+\de)}$, a cubic critical point $p$ of index $i$
of $f_{t_0}$ and cylinders $\left(B_t\right)$ which are neighborhoods of $p$. They are endowed
with Whitney coordinates $(x,y,z)\in \R\times \R^i\times 
\R^{n-i-1}$ so that 
$f_t\vert B_t$ reads:
$$f_t\vert B_t= x^3- (t-t_0)x-\vert y\vert^2+\vert z\vert^2 + cst.
$$
If $\left(X_t\right)_{t\in (t_0-\de,t_0+\de)}$ is a path of pseudo-gradients in the sense of  
\ref{pseudo},  $X_t\vert B_t$  is required to be the descending gradient  of $f_t$
with respect to the Euclidean metric 
of  the Whitney coordinates for every $t\in (t_0-\de,t_0+\de)$ (not only for $t=t_0$). 

The stable/unstable manifold  $W^{u/s}(p,X_{t_0})$ is described now.
One checks that the $x$-axis is the kernel of the Hessian of $f_{t_0}$. 
The half space $\{(x,y,z)\mid x\leq 0,z=0\}$
is the (local) unstable manifold $W^u(p)$; its boundary is the so-called {\it strong-unstable} manifold. 
Similarly, the  half space $\{(x,y,z)\mid x\geq 0, y=0\}$ is the (local) stable manifold and its boundary is the 
{\it strong-stable} manifold. 

Generically, 
$X_{t_0}$ 
has no $j/i$ connections where $j\leq i$, except for possible $i/i$ connections from $p$ to a critical point of index $i$  at a lower level and these connections do not belong to  the strong-unstable manifold of $p$.  Moreover, the $i+1/i$ connections 
are transverse; so, this will be the case for every $t\in (t_0-\de,t_0+\de)$ if $\de$ is small enough.

Moreover, if $\de$ is small 
with respect to the ``horizontal'' size of the cylinders, the cubic critical point $p$
gives rise to a pair of  Morse critical points $(p_t,q_t)\in B_t$ for every $t\in (0,\de)$:
the point $p_t$ has index $i+1$ and coordinates $\left(-\sqrt\frac{t-t_0}{3}, 0,0\right)$;
the point $q_t$ has index $i$ and coordinates $\left(\sqrt\frac{t-t_0}{3}, 0,0\right)$. The closure
of $W^u(p_t, X_t)\cap B_t$  reads $\left\{x\leq  x(q_t)
, z=0\right\}$. 
The closure of 
$W^s(q_t,X_t)\cap B_t$ reads $\left \{x\geq x(p_t)
 , y=0\right\}$. 
One sees a unique connecting orbit from $p_t$ to $q_t$ and 
all other  orbits in $W^u(p_t)$ (resp. $W^s(q_t)$) intersect the bottom (resp.  the top) of $B_t$, which lies
 in a level set of $f_t$ according to  Definition \ref{birth}.
}
\end{rien}
\vskip .6cm

 {\begin{center}
 
 \hskip 0cm \includegraphics[scale=0.6]{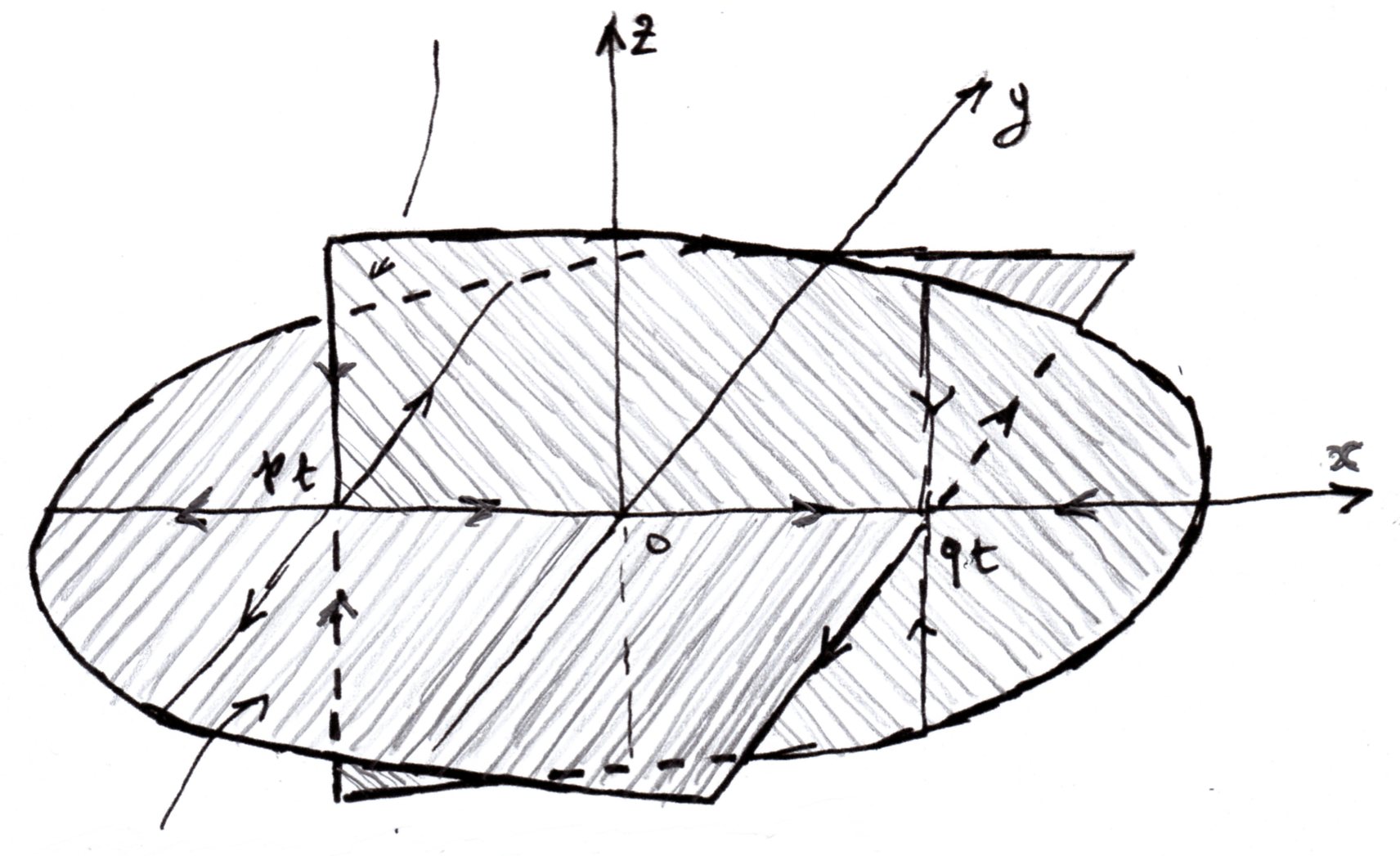}
\end{center}}
\vskip -.4cm ${}$\hskip 2cm {\small $W^u(p_t)\cap \{f\geq f(q_t)-\ep\}$}\\

\vskip -6.6cm{${}$\hskip 2.5 cm {\small $W^s(q_t)\cap\{f\leq f(p_t)+\ep\}$}}

\vskip 6.2cm
\centerline{ Figure 2: After a birth}
\medskip

 \begin{lemme}{\bf (Elementary swallow tail lemma\footnote{In Cerf \cite{cerf} the swallow tail lemma requires no assumption about pseudo-gradient lines but there are some topological assumptions.}).}\label{swallow}
  Let $\ga:=\left(f_t\right)_{t\in [0,1]}$ be a generic path of functions on $M$. 
   Assume that its restriction to $t\in [t_0,t_1]$ has a Cerf graphic showing a \emph{swallow tail}
 as in  figure 3A: there are three critical points, $p_t, p'_t$ of index $i+1$ and $q_t$ of index $i$,
  such that the pair $(p_t,q_t)$ is created at time   $t_0$  and the pair $(p'_t, q_t)$ is cancelled at time $t_1$; at some 
    $\tau\in (t_0,t_1)$  the critical values are equal:
  $f_\tau(p_\tau)=f_\tau(p'_\tau)$. 
  Moreover, it is given a generic  family of pseudo-gradients
  $X_t$ for $f_t$  
  satisfying the next conditions for every $t\in [t_0,t_1]$:
  \begin{itemize}
  \item $W^u(p_t)$ (resp. $W^u(p'_t)$) intersects $W^s(q_t)$ transversely along 
a single  orbit $\ell_t$ (resp. $\ell'_t$);
\item every other orbit  in $W^u(p_t)$  and $W^u(p'_t)$  crosses the  level set   $a_t:=f_t(q_t)-\ep$,
for some $\ep>0$.
\end{itemize}
  Then, given $\de>0$,  
  the path $\ga$ can be deformed   
  to a path $\ga'$ whose Cerf graphic is trivial over 
 $[t_0,t_1]$ as in figure 3B, the deformation being stationary on $[0, t_0-\de]\cup [t_1+\de, 1]$. 
 
  \end{lemme}
 \begin{center}
  \parbox[t]{15cm}{
 \hspace{3.5cm}\includegraphics*[scale=0.5]{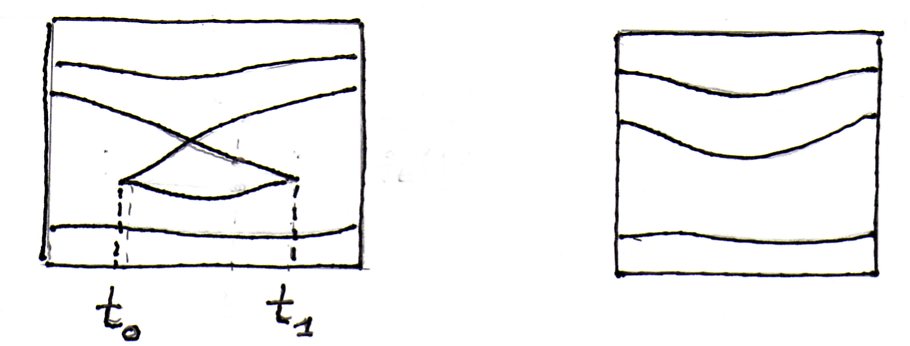}
 
\hspace{4.2cm}Figure 3A\hskip 3cm Figure 3B
}
\end{center}
\vskip .5cm

   \nd {\bf Proof.}  There are three parts.
   
   \nd {\sc A) General setup.}  First, we choose birth cylinders $B_t, \ t\in (t_0-\de', t_0+\de')$ as in
   \ref{fam-pseudo}, the $\de'$ being provisional.
   Without loss of generality, we may assume $f_t\vert B_t= x^3-(t-t_0)x-\vert y\vert^2+\vert z\vert^2$ 
   (no additive constant). 
    And similarly for the cancellation 
   time $t_1$. Take  $\ep$ as in the above statement  and truncate the birth cylinders at level $\pm 2\ep$;
   from  now on, $B_t$  will denote the truncated cylinder.
    
    Set  $\de=\de(\ep)$, 
    so that, for $t=t_0+\de$, the two critical points of  $f_t$  in $B_t$ have value $\pm \ep$.
    Decreasing $\ep$ if necessary, we get  $\de<\de'$. 
   Moreover, 
    except the connecting orbit, every $X_t$-orbit in the invariant manifolds of $p_t$ and
    $q_t$ 
    exits $B_t$ through the top or the bottom of $B_t$. And similarly for the pair 
    $(p'_t, q_t)$ when $t\in [t_1-\de, t_1]$.
    
    Since $f_t(p_t)-f_t(q_t)$ is increasing when $t$ is close to $t_0$,  by taking $\ep$ small enough
     we have
    $f_t(p_t)-f_t(q_t)>2\ep$ for every $t\in (t_0+\de, t_1]$. Similarly, $f_t(p'_t)-f_t(q_t)>2\ep$ 
   for every $t\in [t_0, t_1-\de)$.
    
    For $t\in [t_0+\de,t_1-\de]$, we are going to choose Morse models $\M(q_t), \M(p_t),\M(p'_t)$
    with coordinates $(x,y,z)\in \R\times \R^i\times\R^{n-1-i}$ so that:
    $$\begin{array}{ll}
    f_t\vert\M(q_t)= +x^2-\vert y\vert^2+\vert z\vert^2+f_t(q_t),&\M(q_t)
    \subset f_t^{-1}\left([f_t(q_t)-\ep,f_t(q_t)+\ep]\right)\\
    f_t\vert\M(p_t)= -x^2-\vert y\vert^2+\vert z\vert^2+f_t(p_t),&\M(p_t)
    \subset f_t^{-1}\left([f_t(p_t)-\ep,f_t(p_t)+\ep]\right)\\
     f_t\vert\M(p'_t)= -x^2-\vert y\vert^2+\vert z\vert^2+f_t(p'_t),&\M(q_t)\subset f_t^{-1}\left([f_t(p'_t)-\ep,f_t(p'_t)+\ep]\right)\,.
    \end{array}
    $$
    The pseudo-gradient $X_t$ will be  tangent to the lateral boundary of these models without specifying more.
    Observe that 
    $\M(q_t)$ and $\M(p_t)$ are disjoint for every $t>t_0+\de$; and similarly for 
    $\M(q_t)$ and $\M(p'_t)$ when $t<t_1-\de$.  
    
   We begin by  fixing $\M(p_t)$ and $\M(q_t)$ when $t = t_0+\de$. We choose  their 
    $(y,z)$-coordinates 
    to be those of $B_t$; only the $x$ coordinate has to be changed to have Morse coordinates. 
    And similarly for  $\M(p'_t)$ and $\M(q_t)$ when $t = t_1-\de$.

    Then, we refer to Corollary \ref{unstable-adapted} for extending 
    the choice of  Morse coordinates about $p_t$ to $t>t_0+\de$ 
    so that they are adapted to $\left(f_t, W^u(p_t)\right)$ for every $t$.
    The same is done for $\M(p'_t)$, $t< t_1-\de$.
    For  $\M(q_t)$, $t\in [t_0+\de,t_1-\de]$, we do almost the same except for two differences: 
    \begin{enumerate}
    \item The Morse coordinates are chosen to be  adapted to the stable manifold $W^s(q_t)$.
    \item Since the coordinates 
    are already fixed for $t=t_0+\de$ and $t=t_1-\de$,  item 2 of Corollary \ref{unstable-adapted}
    has to be used.
    \end{enumerate}
    
    Once this choice is made, nothing prevents us from modifying $X_t$ in each  considered Morse model,
    so that it becomes tangent to the $x$-axis,
    the $y$-space and the $z$-space respectively,  as it is the case 
    in $B_t$ when $t\in [t_0, t_0+\de]$ and  $t\in [t_1-\de, t_1]$. The unstable manifolds of $p_t$ and $p'_t$
    are kept unchanged and also the stable manifold of $q_t$; but  the unstable manifold of  $q_t$
    now satisfies
    $$ {\rm (A1)} \quad W^u(q_t)\cap \M(q_t)= \{x=0, z=0\}.$$

    We now recall the {\it cut-and-paste} construction for vector fields, which is abundantly used
   in \cite{h-cob} without using this name.
    Given a Morse function $f$ and a pseudo-gradient $X$, the change of $X$ by {\it 
    cut-and-paste} along a regular  level set $\{f=c\}$ consists of the following:
     cut $M$ at this level, make an isotopy of the upper part
    $\left(\psi_s\right)$ so that $(\psi_1)_*X$ has the same germ as $X$ along the cut, and finally
    glue  $(\psi_1)_*X$ in the upper part to $X$ in the lower part. The assumption for the germs  
    guaranties the smoothness of the resulting vector field. The same construction works in a family.

    By hypothesis of Lemma \ref{swallow},   
    the trace of $W^u(p'_{t})$ in the top of $B_{t},\ t\in [t_0,t_0+\de]$,  intersects 
    transversely   the trace of  $\overline{W^s(q_{t})}$ in  a single point $m_t$. The latter
     trace is a closed disc bounded
    by the trace of $W^s(p_t)$. Moreover, by the genericity assumption in \ref{fam-pseudo}
     the point $m_t$ lies in the interior of that disc. So, we may apply 
     cut-and-paste in the top of $B_t$ to make the part of $W^u(p'_{t})\cap B_t$  
    lying close to  $\{y=0\}$ to be 
    contained in $\{z=0, x>x(q_t)\}$  
     for every $t\in [t_0,t_0+\de]$; this construction extends easily to 
    $t\in (t_0-\de, t_0+\de]$. 
     And similarly for $W^u(p_{t})$ in $B_{t}$ for $t\in [t_1-\de,t_1+\de)$.
     
     In the same way, when $t\in [t_0+\de,t_1-\de]$, 
     cut-and-paste applied in the top of $\M(q_t)$
      makes the part of $\left(W^u(p_t)
     \cup W^u(p'_t)\right)\cap \M(q_t)$ lying near $\{y=0\}$ to be contained in $\{z=0\}$. So, the connecting 
     orbits  cover the $x$-axis of $\M(q_t)$. 
      As the support of the isotopy 
     is located near the stable manifold of $q_t$, the orbits in the unstable manifolds of $p_t$
     and $p'_t$, apart from the connecting orbits, descend to the level $a_t=f_t(q_t)-\ep$. \\

    \nd {\sc Claim 1.} {\it There exists an arc $A_t$  in $M$ 
    passing through $(p_t,q_t,p'_t)$ (or only one of them when a pair  of critical points has disappeared), 
    depending smoothly on 
    $t\in (t_0-\de,t_1+\de)$ such that
    the Cerf graphic of 
    $t\mapsto f_t\vert  A_t$ shows a one-variable swallow tail}.\\
    
    \nd {\sc Proof.}     Starting from the above situation of invariant manifolds, a new  cut-and-paste
     makes $\ell_t$ (resp. $\ell'_t$) 
     coincide with the $x$-axis near the bottom of $\M(p_t)$ (resp. $\M(p'_t)$)
    when $t\in [t_0+\de, t_1-\de]$. 
    
    When $t\in (t_0-\de, t_0+\de]$, $A_t$ is made of the $x$-axis of $B_t$, a piece of $\ell'_t$ from 
    $B_t$ to $\M(p'_t)$, the $x$-axis of $\M(p'_t)$ and a path descending transversely to the level sets 
    from the latter to the level  $f_t(q_t)-\ep$. 
    A similar construction
     is performed on the other intervals of $t$. \bull

   \nd {\sc B) Proof of the swallow tail lemma in case $i=0$.}
  This is the only case needed for proving Theorem \ref{ord-f}.\\

    \nd{\sc Claim 2.} {\it Set $h_t:= f_t\vert A_t$. There are coordinates 
$(x, z)\in \R\times\R^{n-1}$ on a  neighborhood $N_t$ of $A_t$, depending smoothly on 
$t\in (t_0-\de,t_1+\de)$,  
such that  $$\begin{array}{rl}
{(i)}&\quad A_t= \{z=0\} \\
{(ii)} &\quad f_t(x,z)=h_t(x)+\vert z\vert^2.
\end{array}
$$
}
\nd{\sc Proof.} Indeed, it is true on a neighborhood $U_t$ of the set of critical points $\{p_t,p'_t,q_t\}$ 
by the choice we made of the Morse models in A).
First, extend this coordinates arbitrarily
so that $(i)$ holds.
As $h_t$ restricted to 
 $A_t\smallsetminus U_t$ has no critical points, Lemma \ref{conjugation} 
 applies with one parameter 
 $t\in (t_0-\de,t_1+\de)$ and the 
 following correspondence of notation: $V=M$, 
 $V'= A_t\smallsetminus U_t$,  $K=\partial V'$, $f=f_t$, $g= h_t+\vert\cdot\vert^2$.
 ${}$\bull\\

 Now, choose a function $h^1_t$ coinciding with $h_t$ near 
 the boundary of $A_t$  with a single  critical point, indeed a maximum, and satisfying $h^1_t(x)\leq h_t(x)$ for every $x\in A_t$.
   For $s\in [0,1]$, set $k_t^s(x)=s\left(h^1_t(x)-h_t(x)\right)$ and consider the deformation of path of functions $s\mapsto \left(h_t^s\right)_t$ given by 
   $$ (*)\quad \quad h_{t}^s(x)=h_t(x)+k_t^s(x).
   $$
   Note that the path $\left(h_t^1\right)$ has a ``trivial''
   Cerf graphic. So,  the formula $(*)$ solves  the {\it one-dimensional} 
   elementary swallow tail lemma.
   
   Using the coordinates given by  Claim 2, the deformation  extends to the neighborhoods $N_t$ thanks to the formula
  $$s\mapsto h_t(x)+ \om(\vert z\vert )k^s_{t}(x) +\vert z\vert^2, 
$$  where $\om$ is a bump function with a small support,
  centered at 0. The $z$-derivative vanishes at 
  $z=0$ only and the critical points are those of the one-dimensional case. Moreover,
  the deformation is stationary on the boundary of $N_t$ and, hence, extends to $M$
  as  a family $s\mapsto \left(f_t^s\right)_{t\in(t_0-\de,t_1+\de)}$.
  When 
  $s=1$, the Cerf graphic of $\left(f_t^s\right)_{t\in[t_0-\de,t_1+\de]}$ is trivial and the swallow tail
  lemma is proved when $i=0$. 
  ${}$ \bull
  
  \nd {\sc C) Proof of Lemma of the swallow tail lemma when $i>0$.} 
  We continue with the birth cylinders and 
  the Morse models we introduced  in part A).\\
  
  \nd{\sc Claim 3.} {\it There exists a smooth one-parameter family $\left(W_t\right)_{t\in(t_0-\de,t_1+\de)}$
  of smooth compact $(i+1)$-submanifolds, such that:
  \begin{itemize}
  \item $A_t\subset W_t$,
  \item $\partial W_t$ lies at  
   level 
    $a_t$  of the end points of $A_t$,
    \item the only critical points of $f_t\vert W_t$ are $p_t,q_t, p'_t$ and are non-degenerate except for 
    the cubic points when $t$ equals $t_0$ or $t_1$.\\
    \end{itemize}
  }
  
  \nd{\sc Proof.} As a consequence of the cut-and-paste we have made, 
  the closure of $W^u(p_t)$ in the upper 
  level set $\{f_t\geq a_t\}$ and the one of $W^u(p'_t)$ 
  intersect precisely  
  the part of $W^u(q_t)$ lying in that upper level set. Moreover, both match smoothly along
  this common part of  their 
  boundary. This is given for free by the last choice of  pseudo-gradients (see Formula (A1)). 
  So, we set  
  $$W_t= \left[W^u(p_t)\cup W^u(q_t)\cup W^u(p'_t)\right]\cap \{f_t\geq a_t\}.
  $$ 
\bull
   
   \nd{\sc Claim 4.} {\it There are coordinates 
$(x, y,z)\in \R\times\R^i\times \R^{n-i-1}$ on a  neighborhood $N_t$ of $A_t$, depending smoothly on 
$t\in (t_0-\de,t_1+\de)$,  
such that  $$\begin{array}{rl}
{(i)}&\quad A_t= \{y=0,z=0\} \ {\rm and}\  W_t= \{z=0\},\\
{(ii)} &\quad f_t(x,y,z)=h_t(x)-\vert y\vert^2+\vert z\vert^2.
\end{array}
$$
}\\

\nd{\sc Proof.} This is similar to  Claim 2, except that here Lemma \ref{conjugation} has to be applied twice:
firstly in a neighborhood $\mathcal V_t$ of $A_t$ in $W_t$ and secondly in a neighborhood of  
$\mathcal V_t$ in $M$.\bull

    The radial  vector field $Y_t:= \sum_1^i y_j\partial_{y_j}$ in $N_t$ is transverse to the level sets of $f_t$
  in $(N_t\smallsetminus  A_t) \cap\{z=0\}$. 
 Keeping its notation, it extends to  $W_t$ as a Lyapunov vector field (meaning that 
 the Lyapunov inequality holds)
 for $f_t\vert (W_t\smallsetminus A_t)$ since $f_t$ has no critical points on 
 $W_t\smallsetminus  A_t$\,. So, by following the trajectories of $-Y_t$ we get a fibration 
 of $W_t$ over $ A_t$ in $i$-discs, pinched at the end points of $ A_t$ (the diameter of the fibre vanishes there). The fibre
 $D_x$ over $x\in  A_t$ is equipped with a Morse function, namely 
 $g_{t,x}:= f_t\vert D_x$, which has one critical point, a maximum
 indeed, at  $x\in  A_t$.

 Extend $Y_t$  to some neighborhood 
 $\widetilde N_t$ of $W_t$ in $M$ 
 as a Lyapunov vector field $\widetilde Y_t$ of $f_t\vert (\widetilde N_t \smallsetminus A_t)$. 
 Choosing $\widetilde N_t$ to be invariant by the positive semi-flow of $\widetilde Y_t$ gives 
 $\widetilde N_t$ a structure
 of bundle over $A_t$ whose fibre $\widetilde D_x $, $x\in  A_t$,
 is diffeomorphic to $D_x\times D^{n-i-1}$.
 The restriction $\tilde g_{t,x}$ of $f_t$ to the
 fibre $\widetilde D_x$, $x\in  A_t$,
  is a Morse function with the single critical point $x\in A_t$.
  It is  equipped with the pseudo-gradient $\widetilde Y_t$,
 whose  unstable manifold is $D_x$. 
 
 We apply Lemma \ref{reorder} to the function $\tilde g_{t,x}$, where $(t,x)$ is a parameter.
 This lemma allows us to decrease the critical value $f_t(x)$ as we want, 
  without introducing new critical points, as long as this value remains greater than 
  $f_t(\partial W_t)= a_t$. This process yields a deformation of $\left(f_t\right)$
  which extends the solution $(*)$ of the one-dimensional swallow tail lemma
  without introducing new critical points, and solves the general case. \bull\\

   \begin{lemme}{\bf (Elementary lips lemma).}\label{lips}
   Let $\ga:=\left(f_t\right)_{t\in [0,1]}$ be a generic path of functions on the manifold $M$.
Assume that its restriction to $t\in [t_0,t_1]$ has a Cerf graphic 
 as in figure 4 \emph{(lips)}:  for $t\in (t_0,t_1)$, there are two critical points $p_t, q_t$ of respective indices
 $i+1$ and $i$
such that the pair $(p_t,q_t)$ is created at time   $t_0$  and is cancelled at time $t_1$.
Moreover, a smooth family of pseudo-gradients $X_t$ for $f_t$ is given satisfying the next
conditions for all $t\in [t_0,t_1]$:
  \begin{itemize}
  \item $W^u(p_t)$ intersects $W^s(q_t)$ transversely   along 
a single orbit $\ell_t$;
\item all the other orbits  in $W^u(p_t)$  cross the  level set   $f(q_t)-\ep$,
for some $\ep>0$.
\end{itemize}
   Then $\ga$ can be deformed  
  to a path $\ga'$ so that the corresponding  lips are removed from the  Cerf graphic, the deformation being stationary on $[0, t_0-\de]\cup [t_1+\de, 1]$
 for any $\de>0$. \\
  \end{lemme} 
   \begin{center}
  \parbox[t]{15cm}{
\hspace{4cm} \includegraphics*[scale=0.5]{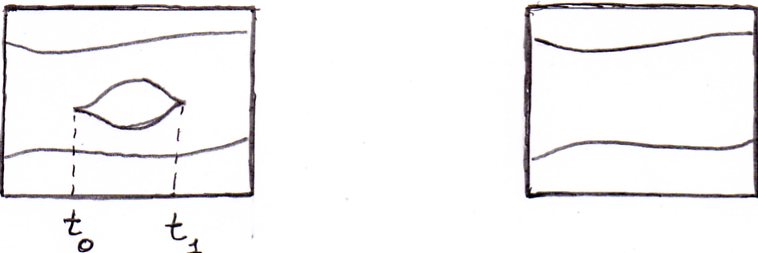}

\centerline{Figure 4A \hskip 2.5cm Figure 4B}
}
\end{center}
\vskip 1cm

   \begin{lemme}{\bf (Morse's cancellation theorem).}\label{cancel}
   Let $f:M\to \R$ be a Morse function equipped with a pseudo-gradient $X$.
   Let $(p,q)$  be a pair of critical points of consecutive indices whose invariant manifolds
   satisfy the next conditions:
   \begin{itemize}
   \item $W^u(p)$ intersects $W^s(q)$ transversely and along a single  orbit ;
   \item all the other orbits in $W^u(p)$ cross the level set $f(q)-\ep$ for some $\ep>0$.
   \end{itemize}
   Then, for every 
   small neighborhood $U$ of the closure of  the intersection
    $W^u(p)\cap \left\{f\geq f(q)-\ep\right\},$ 
   there is 
   a Morse function which has no critical points in $U$ and  coincides with $f$ away from $U$.\\
   
   \end{lemme} 

  \section{Path of polar functions}\label{s4} 
 
 \begin{rien}{\bf Proof of Theorem \ref{ord-f}.} {\rm According to Theorem \ref{order-path}, there is a
 path $\ga:= \left(f_t\right)$  fulfilling all requirements of Theorem \ref{ord-f} (birth times before cancellation times and order of critical values) except the one min/one max condition. So, the matter is to kill the appearance of extra local minima or maxima.
 We are  looking at the local minima only. 
 
First, we make the assumption (H) that  one can follow continuously a minimum $m_t$ of $f_t$ from $t=0$
 to $t=1$. 
 By permuting the birth times if necessary (since $\dim M>1$,  the last claim  of Lemma \ref{shift} applies) and 
 cancelling by pairs the crossings of index 0 critical 
 values (Lemma \ref{reorder}), we may assume that the index 0 part of the Cerf
 graphic shows no crossings (see figure 5A). 
 
 Let $\mu$ be the 
 maximal number of extra minima along $\ga$; we are going to decrease $\mu$ by 1. 
 Denote $(t'_0,t'_1)$ the interval  where $f_t$ has $\mu$ extra minima.
 For $t\in (t'_0,t'_1)$, denote  the upper local minimum of $f_t$ by $m'_t$.
 
 Without loss of generality we may assume that 3/2 separates the index 1 critical values from those of index 2; the same is true for the value $3/2-\eta$, if $\eta>0$ is small. 
 Set
 $L_t:=f^{-1}_t(3/2-\eta)$.
  Since $M$ is connected and  $L_t$ lies above all the critical points of index 1, $L_t$
 is connected.
 
  If $X_t$
 is a pseudo-gradient of $f_t$, we see  in $L_t$ the trace $S_t$ of the stable manifold
 $W^s(m_t,X_t)$ and, when $t\in (t'_0,t'_1)$, the trace $S'_t$ of the stable manifold
 $W^s (m'_t,X_t)$. Both are changing when  handle slides of index 1 happen. But, due to 
 $n\geq 3$, they remain connected; indeed,  each one is always an $(n-1)$-sphere with holes. 
 
  \begin{center}
  \parbox[t]{15cm}{
 \hspace{4cm}\includegraphics*[scale=0.5]{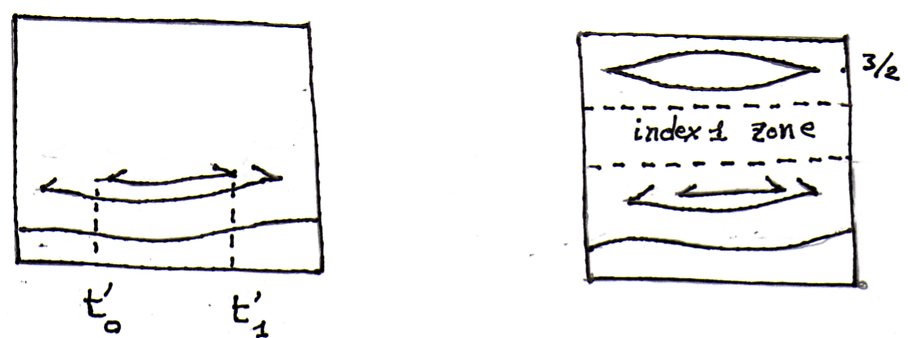}
 
 \vspace{0cm}
\centerline{\rm Figure 5A\hskip 3cm Figure 5B}
}
\end{center}

 
 So, choose smoothly 
 points $x_t\in S_t$ 
 and $x'_t\in S'_t$ linked by a simple arc $\al_t$ in 
 $L_t$. We introduce a cancelling pair of critical points $(s_t,r_t) $ of respective index $(2, 1)$
  in a collar neigborhood
 above
 $L_t$; the birth time is  chosen less than $t'_0$,  the cancellation time greater than $t'_1$ (compare
 figure  5B),  and the base of the birth cylinder is a $(n-1)$-disc in $L_t$ centered at $x_t$.
   Denote by $\ga' :=\left(f'_t\right)$ this new path from $f_0$ to $f_1$.
 After choosing a suitable pseudo-gradient $X'_t$, we have for every $t\in [t'_0+\ep,t'_1-\ep]$:
 $$W^u(r_t,X'_t) \cap L_t= \{x_t,x'_t\}, \  W^u(s_t, X'_t)\cap L_t=\al_t\,.
 $$
 In particular, there are no $X'_t$-connecting orbits form $r_t$ to another   critical point of index 1. Therefore, Lemma \ref{reorder} applies and a new deformation
 of the path $\ga'$ puts the critical value of $r_t $ below the other critical values of index 1
 when $t\in [t'_0+2\ep,t'_1-2\ep]$
 (compare  the Cerf graphic in figure 6A).
 By the choice of $x'_t$, 
 there is exactly one connecting orbit from $r_t$ to $m'_t$ for every $t\in [t'_0+2\ep,t'_1-2\ep]$. 
 One makes  
 cancellations at times $t'_0+2\ep$ and $t'_1-2\ep$.
 These cancellations may be viewed as a new deformation of the path $\gamma'$;
  the final Cerf graphic looks like figure 6B, with two swallow tails separated by lips. Lemma \ref{swallow} and \ref{lips} apply
   and yield  some deformation of the path  of functions so that 
   the swallow tails  and lips vanish. The final path of this last deformation 
   has $\mu-1$ extra minima. This finishes the proof in case of (H).
}\end{rien}

\begin{center}
\parbox[t]{15cm}{
\hspace{3.5cm} {\includegraphics*[scale=0.4]{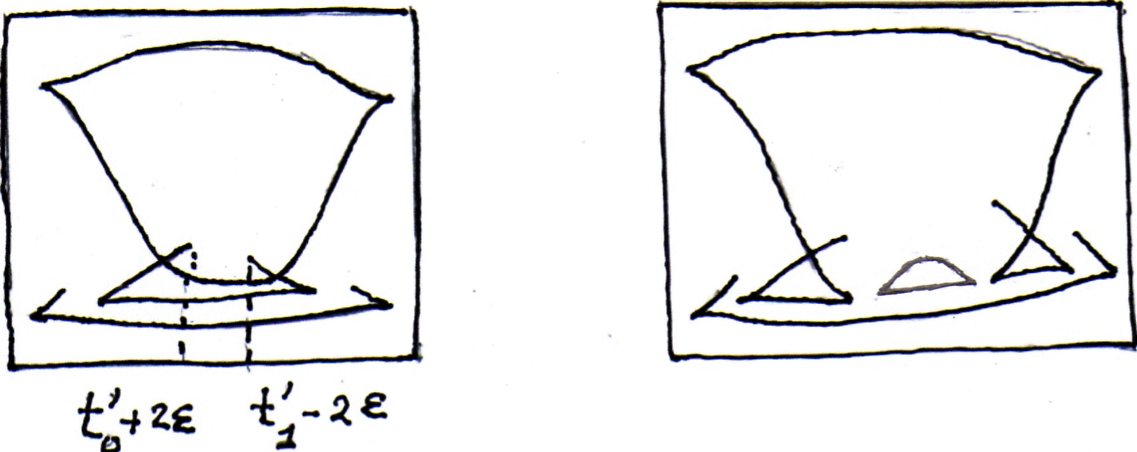}}

\centerline{\rm Figure 6A\hskip 3cm Figure 6B}
}
\end{center}
 \vskip 1cm
 
I am indebted to the anonymous referee who made me observe that the general case easily reduces to assumption (H). Indeed, a suitable isotopy of $M$ makes the minima (resp.  
maxima) of $f_0$ and $f_1$ coincide. Since the germ of smooth
function is unique at a non-degenerate extremum, up to isotopy and 
rescaling,  we may assume that $f_0$ and $f_1$ coincide 
on small discs $d$ and $d'$ about these extrema. Then, by connecting $f_0$ to $f_1$ in the space
of smooth functions having a given restriction to $d$ and $d'$, (H) is fulfilled.\bull
 
\vskip 1cm
\begin{rien} {\bf Final comments.}{\rm

\nd 1) The Reidemeister-Singer theorem, that is, item 1 in Theorem \ref{r-s}, is also proved by 
 R. Craggs  in the piecewise linear category (see \cite{craggs}). His proof relies of previous results on collapsings, due to Chillingworth \cite{chilling}. But the original proof was revisited and explained by 
 L. Siebenmann in \cite{siebenmann}.


\nd 2) It is worth noticing that both parts of Theorem \ref{r-s} are consequence of two statements
(Theorems \ref{order-path} and \ref{ord-f})
about functions which hold true in any dimension. 
These two theorems should be known to specialists. Maybe, the proof of Theorem \ref{order-path} that is given here
is almost the simplest one.   I did not find any written proof of Theorem \ref{ord-f}.

\nd 3) The  proof of the latter theorem is not very elementary, due to the use
  of the swallow tail 
lemma. 
So,  the classical 3-dimensional proof of item 2 in Theorem \ref{r-s} remains competitive. The statement reads as this:} Let $H$ be a 3-dimensional handlebody of genus $g$, and let
$\mathcal D, \mathcal D'$ be two  minimal systems of $g$ compression discs of $H$ whose complement
 is a 3- ball.
 Then, one can pass  from $\mathcal D$ to $\mathcal D'$ by finitely many
handle slides. {\rm This can be proved by a very standard} cut-and-past {\rm technique}.
\end{rien}

I am grateful to Francis Bonahon,  Jean Cerf, Alexis Marin and Patrick Massot 
for comments on  first versions of this note. I am indebted to the referee who suggested me 
several improvements. Carlos Moraga Ferr\'andiz \cite{carlos} is the first who  used of the techniques introduced in this note; I thank him for valuable suggestions. I am also grateful to Marc Chaperon for 
discussions about the saddle-node bifurcation.\\

\end{document}